\documentclass[12pt]{article}

\usepackage{amssymb}
\usepackage{amsthm}
\usepackage{amsmath}
\usepackage{authblk}  
\usepackage[margin=3.2cm]{geometry}

\newtheorem{theorem}{Theorem}

\newcommand{\boldx}{\mathbf{x}}

\title{Optimizing feedstock imports with environmental constraints}

\date{}

\author[1]{Massimo Frittelli}
\author[2]{Pierluigi Toma}

\affil[1]{Department of Engineering for Innovation, \hspace{60mm} University of Salento, 73100 Lecce (Italy) \texttt{massimo.frittelli@unisalento.it}}
\affil[2]{Department of Economic Sciences, \hspace{60mm} University of Salento, 73100 Lecce (Italy) \texttt{pierluigi.toma@unisalento.it}}

\begin{document}

\maketitle

\section*{Abstract}
In addressing the problem of commodity production out of feedstock imports, an eco-environmentally rational agent aims at minimizing the cost of feedstock imports and their transportation, but also the water footprint of the feedstock production process and the water scarcity in the exporting countries. The problem is formulated as a nonlinear program. This study proves the existence of solutions and quantitatively demonstrates that transportation costs and non-uniform feedstock characteristics inhibit feedstock interchangeability. Moreover, it is shown that the interplay between water footprint and water scarcity across countries can inhibit or foster feedstock interchangeability.

\section*{Keywords}
OR in the environment and climate change - Water footprint - Transport cost - Nonlinear programming - Elasticity

\section{Introduction}
It is estimated that almost one-fifth of the world market trades in natural resources. Out of this sheer volume, two-thirds deal in fuel and forms the largest share for the globe to run its industrial complex and is merchandized nearly across all the countries. Marketing fuel is mainly substantial for those countries whose economy runs on its sale worldwide. Twenty-one such countries have a homogenous economic apparatus with more than 80\% exports dealing in natural resources. Nine of these countries commit the lion's share of more than 50\% GDP to the export of natural resources. Correspondingly, it is also important for those countries that highly depend on such imports for their economies and lack the means or power to exploit natural resources in their domains. Geography plays a vital role in this case where natural resources are the assets traded widely in the international market. Dealing globally, natural resources have several features that distinguish them from other commodities, and they have a part to play in the policymaking of the related sector. Unequal placement and distribution of natural resources, geographically speaking, have made some countries to dominate the market, while others have to depend on these countries. Approximately 90\% of the world's proven fossil fuel reserves are found in only fifteen countries. Since drawing on such resources cannot be done elsewhere. Hence incentives to formulate policies for the relocation of production are essentially nonexistent \cite{imf2007guide}. Discrepancies of the likes mentioned above have a major impact in shedding light on the international trade \cite{leamer1984sources, lederman2009commodity}. According to old-school trade theories, the economy thrives when trade allows the transfer of resources from regions of surplus supply to regions with excessive demand. However, such static effects on the economy require careful assessment against the more dynamic effects of trade on the consumption of resources. A considerable amount of literature can be found on such dynamic effects of renewable resources. Numerous studies have shown that an open-door policy of exploitation of resources with weak claims over them may have a disastrous impact leading to their depletion. \cite{chichilnisky1994social, brander1997international, brander1998open, karp2001common}. In contrast, the literature is split on the matter concerning the trade of non-renewable resources. The authors of \cite{kemp1984efficiency} have briefly summarized their findings and applied the Heckscher-Ohlin theory to observe the sustainability in a setting of Hotelling \cite{hotelling1931economics} in a case where producers take into consideration the opportunity cost of a resource that is on the brink of depletion. On the other hand, one more literature search of several abstracts found elements of international trade that concentrate closely on the ultimate route taken by exporters in an imbalanced competition.
The era between 1900 and 1955 saw a dramatic surge in the share of natural resources in international trade, followed by a dip for some decades and then an increase once again. Several factors have been attributed to the long period of this development. Predominantly, these include population increase, industrialization, and declining costs of shipping and transportation. The last decade has seen an escalation in taxes in order to control pollution, which has affected global trade as well. In contrast to this fast expansion of literature, little attention has been paid to control and impact of trade of renewable resources \cite{bulte2005trade}. Failure to adapt to external influence and weak or faulty right claims over the resources means that their management characteristically takes place only in a ``second-best world''. Thanks to the genius of Lipsey \& Lancaster \cite{lipsey1956general}, economists in the present era world over recognize that liberalization in trade in the setting of predetermined contortions may produce uncertain outcomes. The ``second-best'' character of resource management distinctly makes it a topic of interest for economic research on the impact of trade liberalization.\\
The purpose of the present work is to provide a mathematical model that explains the complex process that drives the strategical choices in commodity productions, including economic, environmental, and ethical factors.\\
Relying on underlying theoretical assumptions, the proposed model takes the form of a nonlinear optimization problem. Thanks to the presence of nonlinearities, the proposed model exhibits productive behaviour that reflects the complicated relationship between the heterogeneous factors the model accounts for.\\
The problem of feedstock optimization via programming problems has been addressed primarily for biofuels \cite{miglietta2018optimization, shastri2011development, arifeen2007optimization}, but also for other commodity classes \cite{zhang, marchetti2013multiproduct}. 
In the present work, we propose a general model that is not limited to a single commodity class but aims at grasping a typical underlying behaviour or feedstock optimization in commodity production. In this regard, the present work serves as a theoretical explanatory tool. On the other hand, the proposed model is conceived as an in-depth extension of the linear optimization model for biodiesel feedstock imports previously introduced in \cite{miglietta2018optimization}, to adapt to more general and realistic scenarios. Hence, a possible by-product of the proposed work is an improved decision-making tool for biodiesel feedstock optimization.\\
By carrying out a mathematical analysis of the proposed model, it is concluded that (i) validate the effectiveness of the model and (ii) translate into policy implications.

\section{Problem statement}
\label{}

It is assumed that the agent must produce a fixed quantity of $Q$ of a given commodity by using several interchangeable raw materials imported from several countries.
Specifically, there are $N \in \mathbb{N}$ feedstock-country combinations. For each $i=1,\dots,N$, let $x_i$ be the amount of the $i-th$ country-feedstock that the agent can possibly import and let $\lambda_i$ be the coefficient of transformation of the $i$-th feedstock to the finished commodity. The following constraint is obtained:
\begin{equation}
\sum_{i=1}^N \lambda_i x_i = Q.
\end{equation}
As for the cost, it is assumed that each feedstock $i=1,\dots,N$ has a linear cost $c_i x_i$, plus a sublinear transportation cost $C_i x_i^\gamma$, with $0<\gamma<1$. Hence, the agent must face a total cost of
\begin{equation}
\label{total_cost}
\sum_{i=1}^N (c_i x_i + C_i x_i^\gamma).
\end{equation}
Finally, the impact of water consumption is quantified. Each feedstock $i=1,\dots,N$ requires a virtual water $\mu_i x_i$, where $\mu_i$ is the water footprint (WF) \cite{hoekstra2011water}. Hence, $\mu_i x_i$ could be used as a water consumption indicator. However, in countries affected by water scarcity, the withdrawal of $\mu_i x_i$ liters of water might significantly shrink the natural water reservoir of that country. For this reason, the water consumption is weighted by the normalized remainder of water available in that country, as in the following water impact function:
\begin{equation}
\label{total_water_consumption}
\sum_{i=1}^N \frac{W_i}{W_i - \mu_i x_i} \mu_i x_i,
\end{equation}
where $W_i$ is an estimate of the overall water available in the $i$-th country. Since the agent aims at minimizing the cost and the water consumption impact at once, we are facing a \emph{multicriterial} optimization problem. In order to combine the total cost \eqref{total_cost} and the water impact \eqref{total_water_consumption} into a unique cost function, constant elasticity of substitution (CES) is utilized as aggregator function defined by $CES(a,b) := (a^r + b^r)^\frac{1}{r}$, with $r > 0$, see \cite{heathfield2016introduction}. Hence, the function to minimize is
\begin{equation}
F(x_1,\dots, x_N) := \left(\left(\sum_{i=1}^N  (c_i x_i + C_i x_i^\gamma)\right)^r  + \left(\sum_{i=1}^N \frac{W_i}{W_i - \mu_i x_i}\mu_i x_i\right)^r\right)^\frac{1}{r}.
\end{equation}
The following nonlinear program is obtained:
\begin{align}
\label{nonlinear_program_min_f}
\min\ &  F(x_1,\dots, x_N) \quad \text{s.t.} \\
\label{nonlinear_program_constraint}
&\sum_{i=1}^N  \lambda_i x_i = Q\\
\label{nonlinear_program_positivity}
&x_i \geq 0 \quad i = 1,\dots, N.
\end{align}
This nonlinear program generalizes the linear optimization problem considered in \cite{miglietta2018optimization}, in which transport costs and water scarcity were neglected.

\section{Results}
This study is initiated by proving a necessary and sufficient condition for the existence of one or more solutions to the program \eqref{nonlinear_program_min_f}-\eqref{nonlinear_program_positivity}. Specifically, in the following theorem, it is proved that one or more solutions (i.e. feedstock mixes that are feasible and optimal) exist if and only if the combined water reservoir of the countries of export is sufficient to produce the feedstock needed in the production of the final commodity. This reasonable existence condition shows the well-posedness of the proposed model.
\begin{theorem}[Existence of solutions]
\label{thm:existence}
The nonlinear program \eqref{nonlinear_program_min_f}-\eqref{nonlinear_program_positivity} has at least a solution $(x_1,\dots, x_N)$ if and only if the parameters fulfil
\begin{equation}
\label{parameters_condition}
\sum_{i=1}^N \lambda_i \frac{W_i}{\mu_i} > Q.
\end{equation}
\begin{proof}
The feasible domain of $F$ is the non-compact set
\begin{equation}
D := \left\{(x_1,\dots,x_N)\ \middle|\  \sum_{i=1}^N  \lambda_i x_i = Q, 0 \leq  x_i < \frac{W_i}{\mu_i} \ \forall i=1,\dots,N\right\} \subset \mathbb{R}^N.
\end{equation}
It is proved that $D$ is non-empty if and only if \eqref{parameters_condition} holds true. In fact, if \eqref{parameters_condition} holds true, there exists $0 < \eta < 1$ such that $\sum_{i=1}^N \lambda_i \eta \frac{W_i}{\mu_i} = Q$, hence the point $\left(\eta\frac{W_1}{\mu_1}, \dots, \eta\frac{W_N}{\mu_N}\right)$ is in $D$ and $D$ is non-empty. Conversely, if \eqref{parameters_condition} does not hold and $(x_1,\dots, x_N)$ is in $D$, then $\sum_{i=1}^N \lambda_i x_i < \sum_{i=1}^N \lambda_i\frac{W_i}{\mu_i} \leq Q$, a contradiction, hence $D$ is empty.\\
Finally, the authors were left to prove that, when $D$ is non-empty, problem \eqref{nonlinear_program_min_f}-\eqref{nonlinear_program_positivity} has at least a solution. For any $\varepsilon > 0$, let $D_\varepsilon$ be the compact set defined by
\begin{equation}
D_\varepsilon := \left\{(x_1,\dots,x_N)\ \middle|\  \sum_{i=1}^N  \lambda_i x_i = Q, 0 \leq x_i \leq \frac{W_i}{\mu_i} - \varepsilon \ \forall i=1,\dots,N\right\} \subseteq D.
\end{equation}
Firstly, since $D$ is non-empty, there exists $\varepsilon_0>0$ such that, for $0 < \varepsilon < \varepsilon_0$, $D_\varepsilon$ is non empty as well. 
Now take $\boldx \in D$ and, if $D \setminus D_\varepsilon$ is non-empty, take $\bar{\boldx} = (\bar{x}_1,\dots,\bar{x}_N)\in D \setminus D_\varepsilon$. If $\varepsilon$ further fulfils $\varepsilon < \frac{W_1^2}{W_1 + F(\boldx)}$, then we have $F(\bar{\boldx}) \geq \frac{W_1}{W_1 - \mu_1 \bar{x}_1}\mu_1 \bar{x}_1 \geq F(\boldx)$. Hence, minimizing $F$ on $D$ is equivalent to minimizing $F$ on $D_\varepsilon$, a problem that has a solution from Weierstrass's theorem, since $F$ is continuous on the non-empty compact set $D_\varepsilon$.
\end{proof}
\end{theorem}

\noindent
The following results are devoted to policy implications. To this end, the notion of \emph{productive potential} is introduced, as a quality indicator of feedstocks, based on the unitary price $c_i$, the water footprint $\mu_i$ and the coefficient of transformation $\lambda_i$. For each feedstock $i=1,\dots, N$, its productive potential $P_i$ is defined as
\begin{equation}
P_i := \frac{c_i+\mu_i}{\lambda_i}.
\end{equation}
In the following theorem, a particular case is considered free of transport, no water scarcity, and linear aggregator function. In this particular case, it is proved that there exists an optimal feedstock mix that involves all the considered feedstocks if and only if these possess the same productive potential.

\begin{theorem}
If
\begin{align}
&r = 1 \qquad \text{(linear CES)}\\
&C_i = 0 \qquad \forall i=1,\dots,N \qquad \text{(free transport)}\\
&W_i \rightarrow +\infty \qquad \forall i=1,\dots,N \qquad \text{(no water scarcity)}
\end{align}
then there exists a solution such that $x_i > 0$ for all $i=1,\dots,N$ (i.e. with all country-feedstock combinations involved) if and only if all feedstocks have the same productive potential, i.e.
\begin{equation}
P_i = P_j, \qquad \forall \ i,j=1,\dots,N,\ i \neq j.
\end{equation}
\begin{proof}
If a solution $(x_1,\dots, x_N)$ fulfils $x_i > 0$ for all $i=1,\dots, N$, then it must be a critical point of the decisional function $F$ on the constraint \eqref{nonlinear_program_constraint}. The LaGrange function on such constraint reads
\begin{equation}
\mathcal{L}(\boldx,\xi) = \sum_{i=1}^N c_i x_i + \sum_{i=1}^N \mu_i x_i - \xi\left(\sum_{i=1}^N \lambda_i x_i - Q\right),
\end{equation}
where $\xi$ is the Lagrange multiplier. Hence, the KKT necessary condition is
\begin{equation}
\label{theorem_1_kuhn_tucker}
c_i + \mu_i -\xi \lambda_i = 0 \qquad \forall i=1,\dots,N,
\end{equation}
that is fulfilled if the number
\begin{equation}
\label{theorem_1_condition}
\xi := \frac{c_i + \mu_i}{\lambda_i}
\end{equation}
does not depend on $i$.
Vice-versa, if \eqref{theorem_1_condition} holds true, then $F$ is constant on the feasible set $D$ since
\begin{equation}
F(x_1,\dots,x_N) = \sum_{i=1}^N (c_i + \mu_i) x_i = \xi \sum_{i=1}^N \lambda_i x_i = \xi Q,
\end{equation}
hence every point of $D$ is a solution. Specifically, every point of $D$ with $x_i>0$ for all $i$ is a solution.
\end{proof}
\end{theorem}

\noindent
In the remainder of this work, the following notations will be used:
\begin{equation}
\label{quantities}
\bar{\xi} := \max_i P_i, \qquad \bar{i} := \arg\max_i P_i.
\end{equation}
The index $\bar{i}$ is the feedstock that possesses the maximum production potential $\bar{\xi}$.
In the next theorem, a more general case is considered by including transport costs. In contrast with the previous case, it is proved that \emph{regardless of the production potentials of the given feedstocks}, an optimal mix that involves all the feedstocks, cannot exist. In other words, in the absence of water scarcity, transport costs prevent feedstock interchangeability \emph{a-priori}; that is, one or more feedstocks are necessarily excluded from the optimal mix(es).

\begin{theorem}
If
\begin{align}
&r = 1 \qquad \text{(linear CES)}\\
&W_i \rightarrow +\infty \qquad \forall i=1,\dots,N \qquad \text{(no water scarcity)}
\end{align}
then there exists no solution such that $x_i > 0$ for all $i=1,\dots,N$ regardless of feedstock properties.
\begin{proof}
Suppose by contradiction that a solution $\boldx = (x_1,\dots, x_N)$ fulfils $x_i > 0$ for all $i=1,\dots, N$, then it must be a critical point of the decisional function $F$ on the constraint \eqref{nonlinear_program_constraint}. The LaGrange function on the constraint \eqref{nonlinear_program_constraint} reads
\begin{equation}
\mathcal{L}(\boldx,\xi) = \sum_{i=1}^N (c_i x_i + C_i x_i^\gamma) + \sum_{i=1}^N \mu_i x_i - \xi\left(\sum_{i=1}^N \lambda_i x_i - Q\right),
\end{equation}
where $\xi$ is the Lagrange multiplier. Hence, the KKT necessary condition is
\begin{equation}
\label{theorem_2_kuhn_tucker}
c_i + \gamma C_i x_i^{\gamma-1} + \mu_i  -\xi \lambda_i = 0 \qquad \forall i=1,\dots,N.
\end{equation}
By solving \eqref{theorem_2_kuhn_tucker} for $x_i = x_i(\xi)$, we have
\begin{equation}
x_i(\xi) = \left( \frac{\xi \lambda_i - c_i - \mu_i}{\gamma C_i} \right)^{\frac{1}{\gamma-1}} \qquad \forall i=1,\dots,N.
\end{equation}
The quantity between brackets must be strictly positive, which implies that $\xi > \bar{\xi}$. Observe that the function
\begin{equation}
P(\xi) := \sum_{i=1}^N \lambda_i x_i(\xi) = \sum_{i=1}^N \lambda_i \left( \frac{\xi \lambda_i - c_i - \mu_i}{\gamma C_i} \right)^{\frac{1}{\gamma-1}}
\end{equation}
is well-posed for $\xi \in\ ] \bar{\xi}, +\infty[$, is strictly decreasing w.r.t. $\xi$, tends to $+\infty$ for $\xi \rightarrow \bar{\xi}$ and tends to $0$ for $\xi \rightarrow + \infty$. Hence, there exists a unique $\xi > \bar{\xi}$ such that $P(\xi) = Q$. Since $x_i(\xi) > 0$ for $\xi > \bar{\xi}$, it was found that there exists a unique critical point $\boldx$ of $F$ on the constraint \eqref{nonlinear_program_constraint}. This point is a local maximum since the Hessian of $F$ on $\boldx$ is
\begin{equation}
\mathcal{H}_F(\boldx) = 
\gamma(\gamma -1)
\left[
\begin{matrix}
C_1 x_1^{\gamma-2} & \ & \ \\
\ & \ddots & \ \\
\ & \ & C_N x_N^{\gamma-2}
\end{matrix}
\right],
\end{equation}
which is negative-definite. Hence, $\boldx$ is not a solution of the considered program, a contradiction.
\end{proof}
\end{theorem}

\noindent
In the next and last theorem, the case of free transport is considered but including water scarcity. In this case, the existence of an interplay is showed between water scarcity and water footprint in determining feedstock interchangeability. For each feedstock, reduced water scarcity can compensate for high WF. Specifically, if water scarcity and water footprint fulfill a suitable \emph{compensation condition}, it can be said that:
\begin{itemize}
\item There exists an optimal feedstock mix that involves all the considered feedstocks;
\item The optimal mix is unique. Hence no feedstock can be excluded from the optimal mix.
\end{itemize}

\begin{theorem}
Suppose that
\begin{align}
&r = 1 \qquad \text{(linear CES)}\\
&C_i = 0 \qquad  \forall i=1,\dots,N \qquad \text{(free transport)}.
\end{align}
If the existence condition \eqref{parameters_condition} and the following compensation condition\footnote{Condition \eqref{theorem_3_condition}-\eqref{theorem_3_condition_2} is a bound on the variability of the ratios $\frac{c_i}{\lambda_i}$ and $\frac{\mu_i}{\lambda_i}$ that accounts for water scarcity. In particular, the variability in water scarcity can compensate or worsen the variability in water footprint.} holds true
\begin{equation}
\label{theorem_3_condition}
M_i^2\frac{c_{\bar{i}}}{\lambda_{\bar{i}}} + M_i^2\frac{\mu_{\bar{i}}}{\lambda_{\bar{i}}}< M_i^2 \frac{c_i}{\lambda_i} + \frac{\mu_i}{\lambda_i}\qquad \forall i\neq \bar{i},
\end{equation}
where
\begin{equation}
\label{theorem_3_condition_2}
M_i := \max \left\{0, 1 - \frac{Q\mu_i}{W_i\lambda_i(N-1)} \right\} \qquad \forall i \neq \bar{i},
\end{equation}
then the solution is unique and fulfils $x_i > 0$ for all $i=1,\dots,N$.
\begin{proof}
If a solution $(x_1,\dots, x_N)$ fulfils $x_i > 0$ for all $i=1,\dots, N$, then it must be a critical point of $F$ on the constraint \eqref{nonlinear_program_constraint}. The LaGrange function on such constraint reads
\begin{equation}
\mathcal{L}(\boldx,\xi) = \sum_{i=1}^N c_i x_i + \sum_{i=1}^N \frac{W_i}{W_i - \mu_ix_i}\mu_i x_i - \xi\left(\sum_{i=1}^N \lambda_i x_i - Q\right),
\end{equation}
where $\xi$ is the Lagrange multiplier. The KKT necessary condition is
\begin{equation}
\label{theorem_3_kuhn_tucker}
c_i + \frac{\mu_i W_i^2}{(W_i-\mu_i x_i)^2} - \xi\lambda_i = 0, \qquad \forall i=1,\dots,N.
\end{equation}
By solving \eqref{theorem_3_kuhn_tucker} for $x_i = x_i(\xi)$, we have
\begin{equation}
x_i(\xi) = \frac{W_i}{\mu_i}\left(1- \sqrt{\frac{\mu_i}{\xi\lambda_i - c_i}}\right), \qquad \forall i=1,\dots,N.
\end{equation}
Observe that the function
\begin{equation}
P(\xi) := \sum_{i=1}^N \lambda_i x_i(\xi) = \sum_{i=1}^N \lambda_i \frac{W_i}{\mu_i}\left(1- \sqrt{\frac{\mu_i}{\xi\lambda_i - c_i}}\right)
\end{equation}
is well-posed (and all its summands are non-negative) for $\xi \in [ \bar{\xi}, +\infty[$, is strictly increasing w.r.t. $\xi$, and for $\xi \rightarrow +\infty$ tends to $\sum_{i=1}^N \lambda_i \frac{W_i}{\mu_i} > Q$ from \eqref{parameters_condition}. Now, there exists a unique solution $(x_1(\xi), \dots, x_N(\xi))$ with $x_i(\xi) > 0$ for all $i$ if and only if there exists a unique $\xi > \bar{\xi}$ such that $P(\xi) = Q$, which in turn is true if and only if $P(\bar{\xi}) < Q$.
To this end, observe that the $\bar{i}$-th summand of $P(\bar{\xi})$ is $0$, hence it is sufficient to enforce that each of the other $N-1$ summands does not exceed $\frac{Q}{N-1}$, which yields \eqref{theorem_3_condition}. We have found that there exists a unique critical point $\boldx(\xi)$ of $F$ on the constraint \eqref{nonlinear_program_constraint}. This point is a global minimum since, for any $\boldx\in D$, the Hessian of $F$ on $\boldx$ is
\begin{equation}
\mathcal{H}_F(\boldx) = 
2
\left[
\begin{matrix}
\frac{W_1^2\mu_1^2}{(W_1 - \mu_1 x_1)^3} & \ & \ \\
\ & \ddots & \ \\
\ & \ & \frac{W_N^2\mu_N^2}{(W_N - \mu_N x_N)^3}
\end{matrix}
\right],
\end{equation}
which is positive-definite. 
\end{proof}
\end{theorem}

\section{Conclusions}
In this article, the general predictive and explanatory tool of the feedstock market is proposed on a worldwide scale. The proposed methodology grasps common underlying factors that drive the market of potentially any class of feedstock, thanks to the relatively basic, but also well-justified assumptions our model relies on. On the other hand, the multiple nonlinearities allow the model to capture a plethora of different scenarios and drive strategic choices that would be otherwise inaccessible. For instance, in the previous work \cite{miglietta2018optimization}, the authors have extracted strategic choices on biodiesel feedstock optimization through a particular linear case of the model proposed here, by leaving out the fine-grained complexity of the optimal solutions that we would observe in the much more general current setting. Moreover, several specific studies carried out for different feedstock classes or other similar worldwide marketing problems can now be addressed in the present framework.

The results obtained in this study provide the following policy insights. First, one or more feasible and optimal mixes of feedstock imports exist if and only if the countries of exports possess a combined water reservoir that allows for the production of the final commodity. Second, feedstocks characterized by high WFs and import costs tend to be excluded from the optimal choices. In fact, in the limit case of free transport and no water scarcity, different feedstocks can be included in an optimal mix only if these are precisely homogeneous in terms of \emph{production potential}. Third, transportation costs reduce the possible feedstocks, which can be included in the optimal mix. In fact, in the limited case of no water scarcity, an optimal mix that includes all feedstocks does not exist, regardless of their production potentials. Fourth, water scarcity acts as a balancing factor that can enhance or compensate for heterogeneity in the production potentials of the raw materials. In fact, in the limit case of free transport, there exists a compensation condition in closed form, under which the optimal mix is unique and includes all feedstocks.

\bibliographystyle{plainurl}
\bibliography{bibliography}
\end{document}